\documentclass[12pt]{article}
\title{Un domaine de Fatou-Bieberbach
\\
\`a plusieurs feuillets
\footnote{AMS MSC: 37F25, 37F05}
}
\author{Claudio Meneghini
}
\begin{document}
\maketitle
\bibliographystyle{plain} 
\parindent=8pt
\pagestyle{myheadings}
\sloppy
\markboth
{\large\tt version provixoira\rm------------------------------------aqueste.tex}
{\large\tt version provixoira\rm------------------------------------aqueste.tex}

\font\cir=wncyb10
\def\pe{\cir\hbox{P}}
\def\CIRC{\mathop{\tt o}\limits}
\def\quan{\vrule height6pt width6pt depth0pt}
\def\QUAN{\ \quan}
\def\R{\hbox{\tt R}}

\newtheorem{definition}{Definition}
\newtheorem{defi}[definition]{D\'efinition}
\newtheorem{lemme}[definition]{Lemme}
\newtheorem{proposition}[definition]{Proposition}
\newtheorem{theoreme}[definition]{Th\'eor\`eme}        
\newtheorem{corollaire}[definition]{Corollaire}  
\newtheorem{remarque}[definition]{Remarque}
  
\font\sdopp=msbm10
\def\ERRE {\sdopp {\hbox{R}}}
\def\CI {\sdopp {\hbox{C}}}
\def\DI {\sdopp {\hbox{D}}}
\def\ENNE{\sdopp {\hbox{N}}}
\def\ZETA{\sdopp {\hbox{Z}}}
\def\PI {\sdopp {\hbox{P}}}
\def\M{\hbox{\tt\large M}}
\def\N{\hbox{\tt\large N}} 
\def\P{\hbox{\boldmath{}$P$\unboldmath}} 
\def\f{\hbox{\large\tt f}} 
\def\F{\hbox{\boldmath{}$F$\unboldmath}} 
\def\H{\hbox{\boldmath{}$H$\unboldmath}} 
\def\h{\hbox{\large\tt h}} 
\def\id{\hbox{\boldmath{}$id$\unboldmath}} 
\def\pphi{\hbox{\boldmath{}$\phi$\unboldmath}} 
\def\Ef{\phi}
\def\ppsi{\hbox{\boldmath{}$\psi$\unboldmath}} 
\def\Ppsi{\hbox{\boldmath{}$\Psi$\unboldmath}} 
\def\labelle #1{\label{#1}}
\def\cIrc{  }
{\small {\bf R\'esum\'e}: nous d\'efinissons 
l'id\'ee 
de bassin
d'attraction 
pour 
l'it\'eration 
des germes holomorphes attractifs
de $\CI^N,0$ 
 et
 la notion de
domaine de Riemann-Fatou-Bieberbach.
Ceci est
un domaine de Riemann   $R$ recouvrant
une region  $\Omega\subset\CI^N$,
et \'egalement $\CI^N$
de fa\c con telle que
les 
fibres 
de la projection vers $\CI^N$
soient compatibles avec celles de la projection
vers $\Omega$.
Enfin, \'etant donn\'e un
endomorphisme de $\CI^N$ admettant 
un point fixe r\'epulsif en $0$ (satisfaisant une hypoth\`ese technique suppl\'ementaire),
nous prouvons que le bassin d'attraction 
du
germe inverse admet un recouvrement par un
domaine de Riemann-Fatou-Bieberbach 
et que, en certains cas, ce domaine
est biholomorphe \`a
$\CI^N$.
}
\section{Introduction}
Rappelons
qu'un domaine de Fatou-Bie\-ber\-bach
est un ouvert de $\CI^N$ biholomorphe \`a $\CI^N$; 
le bassin d'attraction $\Omega$ d'un point fixe d'un {automorphisme} $f$ de $\CI^N$ 
a cette propri\'et\'e dans le cas attractif 
(voir \cite{rudin}, ap\-pen\-di\-ce),
et, parfois, dans le cas des applications tangentes 
\`a l'identit\'e (voir \cite{weickert}).
Pour plus d'exemples de tels domaines, voir
\cite{bedfsmill},
\cite{buzzforn},
\cite{fornsib},
\cite{globevnik},
\cite{kimura},
\cite{myrberg},
\cite{sibony},
\cite{stehle}
et
\cite{stensones}.

Nous proposons une g\'en\'eralisation de ce concept,
 consistant en un domaine de Riemann $(\M,\pi)$ recouvrant un ouvert propre $\Omega$
 de $\CI^N$ 
et \'egalement $\CI^N$: il s'agit donc d'une correspondence holomorphe entre 
$\Omega$ et $\CI^N$.
En outre, les fibres
de la projection vers $\CI^N$
seront compatibles avec celles de la projection
vers $\Omega$: cela permet de voir ce type de construction comme
une sorte de 'correspondence biholomorphe'.

Nous montrerons que, \'etant donn\'e 
$h$ un endomorphisme  de $\CI^N$ avec
un point fixe r\'epulsif r\'egulier en $0$ 
(voir d\'ef. \ref{reg}),
le  bassin potentiel d'attraction {\rm (d\'ef. \ref{basattr})} de $0$  
pour l'inverse local $h^{-1}_0$ peut \^etre recouvert par
un domaine de Riemann $(\M,\pi)$ 
qui sera aussi
un domaine de Riemann 
$(\M,\Psi)$ au-dessus de $\CI^N$.
En outre,
$\Psi(x)=\Psi(y)\Rightarrow \pi(x)=\pi(y)$
et, si $0$ est {\it fortement r\'egulier} (voir 
d\'ef. \ref{reg}), $\M$ est biholomorphe \`a $\CI^N$.

\section{Pr\'eliminaires}

Soient $\N$ et $\M$ des 
vari\'et\'es complexes:
 un {\sf \'el\'ement d'application holomorphe}
de
$\N$ dans $\M$ est une paire 
$\left(U,f    \right)$, 
o\`u
$U$ est un ouvert connexe de
$\N$ 
and $f$ une application holomorphe definie
sur $U$ et \`a valeurs dans $\M$.

Deux \'el\'ements
$\left(U,f    \right)$ et $\left(V,g    \right)$ 
sont {\sf reliables} s'il existe une suite
$ \left\{(U_j,f_j)    \right\}_{j=0,....,n}$,
telle que
$\left(U_0,f_0    \right)=\left(U,f    \right)$, $\left(U_n,f_n    \right)=\left(V,g    \right)$
et,
pour chaque 
$j=0,....,n-1$,
$$
\cases
{
U_j\cap U_{j+1}\not= \emptyset,\cr 
f_{j+1}\vert_{U_j\cap U_{j+1}}=f_{j}
\vert_{U_j\cap U_{j+1}}.
}
$$
On dira aussi que $\left(V,g    \right)$
est un {\sf prolongement analytique} de $\left(U,f    \right)$ et vice versa.

Un {\sf domaine de Riemann} au-dessus d'un
ouvert 
$\Omega\subset\N$  est une paire $\left(R,p    \right)$ 
o\`u $R$ est
une vari\'et\'e complexe et 
$p:R\rightarrow\Omega$ 
un biholomorphisme local 
surjectif
(voir aussi \cite{gunros}, p.43).
Une {\sf extension analytique }
d'un 
{\'el\'ement d'application holomorphe}
 consiste en
un domaine de 
 Riemann connexe
$(S,\pi)$
 au-dessus d'un
ouvert 
$\Omega\subset\N$  
tel que  $U\subset \pi(S)$,
en une 
immersion  holomorphe $j\,\colon\, U\rightarrow S$  telle que $\pi\circ j=\id\vert_{U}$ 
et en une application
 holomorphe 
$F\,\colon\, S\rightarrow \M$
telle que $F\circ j=f$.
Un {\sf morphisme} entre deux extensions analytiques
$
\left(S,\pi,j,F\right)
$
et
$\left(T,\varrho,\ell,G\right)
$
du m\^eme 
\'el\'ement $\left(U,f    \right)$ 
est une application holomorphe
 $h\,\colon\, T\rightarrow S$ 
telle que
 $h\circ\ell=j$.
Un morphisme entre deux extensions analytiques
$
\left(S,\pi,j,F\right)
$
et
$\left(T,\varrho,\ell,G\right)
$
du m\^eme \'el\'ement
est une application non constante,
univoquement d\'etermin\'ee -en $j(U)$
et
donc partout en $S$-, par 
$\ell\circ j^{-1}$. En outre,
$\varrho\circ h=\pi$ et $G\circ h=F$
en $j(U)$ donc partout en $S$.

Le seul morphisme entre une extension analytique et lui m\^eme 
est l'i\-den\-ti\-t\'e,  la
composition de deux morphismes est encore un morphisme; si un morphisme admet 
une application holomorphe comme inverse, elle est encore un morphisme: dans ce cas, nous parlons d'un {\sf isomorphisme} de extensions analytiques.

\begin{defi}
Une extension  analytique 
$S$ de l'\'el\'ement $\left(U,f    \right)$ 
est {\sf  maximale} si, pour chaque extension $\widehat S$
de $\left(U,f    \right)$ il existe un morphisme 
$h\,\colon\, S\rightarrow \widehat S$.
\end{defi}

Remarquons que 
deux extensions maximales du m\^eme
element sont forc\'ement isomorphes
et donc 
l'extension analytique maximale
est unique \`a isomorphismes pr\`es.
 
\begin{theoreme}
Tout \'el\'ement $\left(U,f    \right)$ d'application holomorphe
admet une extension analytique maximale
\QUAN
\end{theoreme}

On pourra  
aussi consulter \cite{narasimhan}, chap. 2;6
ou \cite{malgrange} chap. 1(iv).

Le lemme suivante \'etablit une liaison entre les extensions analytiques maximales
de deux \'el\'ements
qui sont 
 inverses l'un de l'autre: il sera 
utilis\'e dans 
le cas fortement r\'egulier
(lemme \ref{dedans} {\tt (iv)}).
\begin{lemme}
\labelle{inverse}
Soient $({\cal U},f  )$
et $({\cal V},g  )$ deux \'el\'ements
d'applications
holomorphes entre ouverts de 
$\CI^N$,
inverses l'un de
de l'autre; soient
 $(   R,\pi,j,\Phi)$ et $(S,\rho,\ell,\Psi)$
leur extensions analytiques maximales:
alors, si
$\, {\cal C}=\hbox{\{\small\tt points critiques de $\Phi$\}}$
et
${\cal D}=\hbox{\{\small\tt points critiques de $\Psi$\}}$,
on a 
$\Phi(R\setminus{\cal C})
=\rho(S\setminus{\cal D})$.
\end{lemme}
{\bf D\'emonstration.}  
{\bf A)} $\Phi(R\setminus{\cal C})
\subset\rho(S\setminus{\cal D})$: 
soit $\xi\in
 R\setminus{\cal C}$ et $\Phi(\xi)=\eta$:
 il existe 
un voisinage ouvert ${\cal U}_1$ de  $\xi$,
ouverts 
${\cal U}_2\subset\pi({\cal U}_1   )$,
${\cal V}_2\subset \Phi({\cal U}_1   )$ 
et une fonction biholomorphe
$g_2:{\cal V}_2\rightarrow{\cal U}_2$
(avec inverse
$f_2:{\cal U}_2\rightarrow{\cal V}_2$
)
tels que
$({\cal U}_2,f_2   )$ et $({\cal U},f  )$
soient reliables ,
aussi que
$({\cal V}_2,g_2   )$ et $({\cal V},g  )$.
Par construction il existe des immersions 
 holomorphes $\widetilde{\hbox{\j}}:{\cal U}_2
\rightarrow R\hbox{ et }
\widetilde\ell:{\cal V}_2\rightarrow S$
telles que
$\pi\circ\widetilde{\hbox{\j}}=\id$ et $\rho\circ\widetilde\ell=\id$.
Soit
${\cal V}_1=\Phi({\cal U}_1)$ et
$$
\Sigma=\{(x,y   )\in{\cal U}_1\times{\cal V}_2: \Phi(x)=y\}
. 
$$
D\'efinissons $J:{\cal V}_2\rightarrow\Sigma$
en posant $ J(v)=(\widetilde{\hbox{\j}}\circ  g_2(v),v )$.
Or $ (\Sigma,pr_2,J,\pi\circ pr_1   )$  est une extension
analytique de $({\cal V}_2,g_2   )$
car
$ 
\pi\circ pr_1\circ J=\pi\circ\widetilde{\hbox{\j}}\circ g_2=g_2.
$
Mais
$({\cal V}_2,g_2   )$ est reliable avec
$({\cal V},g     )$, donc
$ (\Sigma,pr_2,J,\pi\circ pr_1   )$
est une extension analytique de
$({\cal V},g )$.

Gr\^ace \`a la maximalit\'e, cela
 entra\^\i ne qu'il existe une fonction
holomorphe $h:\Sigma\rightarrow S$ 
telle que
$\rho\circ h=pr_2$, donc
$
\eta=pr_2(\xi,\eta   )=\rho\circ h
(\xi,\eta   )\in\rho(S   )
$.

Enfin, par differentation compos\'ee, 
aucun point de $\rho^{-1}(\eta)$ 
ne peut \^etre critique pour
$\Psi$.

{\bf B)} 
$\Phi(R\setminus{\cal C})
\supset\rho(S\setminus{\cal D})$:
soit
$s\in S$ un point r\'egulier de $\Psi$:
il existe un voisinage
$V$ de $s$
ne
contenant que
points r\'eguliers de 
$\Psi$.
\c Ca signifie que, pour chaque
$s^{\prime}\in V$, il existe un 
\'el\'ement d'application holomorphe
$({\cal V}^{\prime},\widetilde g_{s^{\prime}}   )$
(avec ${\rho(s^{\prime})}\in {\cal V}^{\prime}$)
reliable avec
$({\cal V},g   )$;
en outre, il
existe 
une
immersion holomorphe
$\widetilde\ell:{\cal V}^{\prime}\rightarrow V$.
Par {\bf A)} d\'esormais prouv\'e, 
$\Psi(s)\in\pi(R\setminus{\cal C})$,
donc il existe
$p\in R\setminus{\cal C}$
et un voisinage
$W$ de $p$ dans $R\setminus{\cal C}$ 
tels que 
$\pi(p)=\Psi(s)$ et 
$ \pi^{-1}(\widetilde g({\cal V}^{\prime}
   )   )\bigcap W\not=\emptyset$.
Posons
$ W^{\prime}=\pi^{-1}(\widetilde g({\cal V}^{\prime}
   )   )\bigcap W$: on peut supposer, sans perte de 
g\'en\'eralit\'e, que
$\pi$ soit inversible dans $W^{\prime}$: 
alors il existe une immersion holomorphe
$\widetilde{\hbox{\j}}:\widetilde g({\cal V}^{\prime} 
  )\rightarrow W$.
Donc, pour chaque
 $\zeta\in\widetilde{\hbox{\j}}(\widetilde
 g({\cal V}^{\prime}   )   )$,
il
existe $ \eta\in\widetilde\ell({\cal V}^{\prime}   )$
tel que
$ \Phi(\zeta)=\Phi(\widetilde{\hbox{\j}} \circ \widetilde 
g\circ \rho(\eta)  )$.
Or par la
d\'efinition de  extension analytique, on a
$ \Phi\circ\widetilde{\hbox{\j}}\circ\widetilde 
g=\id$, c'est-\`a-dire 
$ \Phi(\zeta)=\rho(\eta)$.
Consid\'erons maintenant la fonction holomorphe 
$\Xi:W\times V\rightarrow\CI^N$ 
d\'efinie en posant $ \Xi(w,v   )=\Phi(w)-\rho(v)$: 
on a 
$
\Xi
\equiv 0
$ 
dans
${\widetilde{\hbox{\j}}(\widetilde g({\cal V}^{\prime}
   )   )\times\widetilde
\ell({\cal V}^{\prime})}$: 
cet ensemble 
est ouvert
dans $W\times V$, donc
$\Xi\equiv 0$ dans $W\times V$.
Cela implique finalement $ \Phi(p)=\rho(s)$,
ce qui conclut
la d\'emonstration.
\QUAN

Les deux \'enonc\'ees suivants correspondent 
aux
lemmes 1 et 3 de l'appendice de \cite{rudin}:
le deuxieme est connu comme th\'eor\`eme de 
Poincar\'e-Dulac.

Soit $G=(g_1...g_N)$ un automorphisme 
triangulaire inf\'erieur polynomial de $\CI^N$:
$$
\cases
{
g_1(z)=c_1 z_1\cr 
g_2(z)=c_2 z_2+h_2(z_1)\cr
...\cr
g_n(z)=c_N z_N+h_N(z_1...z_{N-1}),
}
$$
o\`u les $c_1...c_ N  $ sont des constantes complexes non nulles,  et chaque $h_i   $ 
est une fonction polynomiale de 
$(z_1...z_{i-1})$; soit $\Delta$ le polydisque unit\'e
 en $\CI^N   $.

\begin{lemme}
(a) il existe $M\in\ENNE,\, \beta\in\ERRE$
tels que $deg(G^{k})\leq M$
et $G^{\cIrc k}(\Delta)\subset\beta^k \Delta   $;
(b) si $\vert c_i\vert<1   $ pour tout $i=1...N$,
$\{G^{\cIrc k}   \}$ converge
 uniform\'ement \`a $0$
sur les compacts
de $\CI^N$
et, pour chaque voisinage $V   $  de 
$0   $,
$\bigcup_{k=1}^{\infty}G^{-k}(V)=\CI^N   $.
\QUAN
\labelle{erreuno}
\end{lemme}
\begin{theoreme} {\tt (Poincar\'e-Dulac)}
Soit $V   $ un voisinage de $0   $ in $\CI^N   $,
$F:V\rightarrow \CI^N  $ une application holomorphe avec $F(0)=0   $ et
$F_*(0)$ triangulaire inf\'erieure
;
supposons que toutes les valeurs propres $\lambda_i   $ de $F_*\vert_0:=A   $ satisfassent $\vert\lambda_i\vert<1   $.
Alors il existe:
{\tt (i)}
un automorphisme triangulaire inf\'erieur polynomial $G$ de $\CI^N $ tel que $G(0)=0   $
et $G_*\vert_0=A   $
{\tt (ii)}
des applications polynomiales $T_m:\CI^N\rightarrow\CI^N   $, avec
$T_m(0)=0   $, ${T_m}_*\vert_0=\id   $
telles que
$G^{-1}\,\hbox{\tt o}\, T_m \,\hbox{\tt o}\,
F-T_m=O(\vert z\vert^m)$, $(m=2,3,4....)   $.
\QUAN
\labelle{erretre}
\end{theoreme}

Le lemme suivant d\'ecrit le comportement d'une application holomorphe au voisinage d'un point fixe
attractif:

\begin{lemme}
Soit $V$ un voisinage de $0$ en $\CI^N$, $F:V\rightarrow\CI^N$
une application holomorphe
admettant un point fixe attractif en $0$: alors
il existe $\alpha<1$ et un voisinage ouvert
$\R\subset V$ de $0$
tel que $F^n(\R)\subset \alpha^n \R$.
\labelle{lattes}
\end{lemme}
{\bf D\'emonstration:} 
gr\^ace au lemme de Schur, on peut supposer,
sans perte de g\'en\'eralit\'e, que 
$$
F_*(0)=
\pmatrix
{
\lambda_1&\ldots&\ldots&\ldots&0\cr
a_{21}&\lambda_2&\ldots&\ldots&0\cr
\vdots&\vdots&\ddots &\ &\vdots\cr
\vdots&\vdots&\ &\ddots &\vdots\cr
a_{N1}&\ldots&\ldots&a_{N\, N-1}&
                          \lambda_N\cr
},
$$
o\`u les $\lambda_k$
(avec $\vert\lambda_1\vert
\leq
\vert\lambda_2\vert
\leq
...
\leq
\vert\lambda_N\vert
$)
sont les valeurs propres de $F_*(0)$
et les $a_{jk}$ des constantes complexes.
Soit
$$
E_{\varepsilon}=
\pmatrix
{
\varepsilon^N&0&\ldots&0\cr
0&\varepsilon^{N-1}&\ldots&0\cr
\vdots&\vdots&\ddots&\vdots\cr
0&\ldots&\ldots&\varepsilon\cr
}:
$$
si $\varepsilon$ est suffisamment petit,
il existe $\alpha<1$ tel que 
$$
\Vert E_{\varepsilon}^{-1}F_*(0)E_{\varepsilon}\Vert
(=\Vert E_{\varepsilon}^{-1}F_*E_{\varepsilon}(0)\Vert)
<\alpha;
$$
il existera alors $\varrho>0$ tel que
$E_{\varepsilon}^{-1}\circ F\circ E_{\varepsilon}(B(0,\varrho))\subset B(0,\varrho)$, donc,
si $p\in B(0,\varrho)$, on a 
$\Vert E_{\varepsilon}^{-1}F_*^nE_{\varepsilon}(p)\Vert<\alpha^n$
et 
$$
\Vert E_{\varepsilon}^{-1}F^nE_{\varepsilon}(p)\Vert<\alpha^n
\Vert p\Vert\leq\alpha^n\varrho,
$$
c'est-\`a-dire 
$E_{\varepsilon}^{-1}\circ F^n\circ E_{\varepsilon}(B(0,\varrho))\subset B(0,\alpha^n\varrho)$;
mais alors, on a 
$$
F^n(E_{\varepsilon}(B(0,\varrho)))
\subset
E_{\varepsilon}(B(0,\alpha^n\varrho))
=
\alpha^n E_{\varepsilon}(B(0,\varrho))
 ;
$$
on conclut en posant $\R=E_{\varepsilon}(B(0,\varrho))$.
\QUAN

\section{Le th\'eor\`eme principal}

\begin{defi}
Soit $(\R,F)$
un \'el\'ement d'application holomorphe avec un point fixe attractif en $0$
et tel que $F^n(\R)\subset\alpha^n\R$, pour
un certain
 $0<\alpha<1$;
on dira que 
$p$  
{\bf est dans le bassin potentiel 
d'attraction } de $0$
pour la dynamique 
de $F$
($p\in BPA(F,0)$ dans la suite)
s'il existe
une suite finie
de  points $\{x_{\nu}\}_{\nu=0...N}$ 
et des prolongements analytiques $(V_{\nu},F_{\nu})$ de $F$ tels que
$x_0=p$, $x_{\nu}\in V_{\nu}$, $F_{\nu}(x_{\nu})=x_{\nu+1}$, $F_{\nu}(V_{\nu})\subset V_{\nu+1}$
et 
$\hbox{\tt\large o}_{\nu=0}^{N}\,
F_{\nu}(V_0)\subset\R
$.
\labelle{basattr}
\end{defi}
Soit maintenant $h$ un endomorphisme de $\CI^N$
avec  un point
fixe r\'epulsif en $0$. 

\begin{defi}
 On dira que le point fixe r\'epulsif
 $0$ est
{\sf r\'egulier} s'il admet un voisinage $\R$
sur lequel $h$ est inversible,
$[h\vert_{\R}]^{-k}\subset\alpha^k\R$ pour $0<\alpha<1$ et,
pour chaque $k$, 
$h^{\cIrc k}$ est localement inversible;
 on appellera $\R$ un {\sf voisinage de r\'egularit\'e} de $0$.
Si, de plus, pour chaque $k$, $h^k$ est un rev\^etement topologique 
on dira que le point fixe est {\sf fortement r\'egulier}.
\labelle{reg}
\end{defi}
\begin{lemme}
Soit $0\in\CI^N$ un point fixe r\'epulsif
r\'egulier pour $h$,
$\R$ un {voisinage de r\'egularit\'e} de $0$
et $F:=[h\vert_{\R}]^{-1}$. Alors
$p\in BPA(F,0)$ si et seulement si
 il existe $k\geq 1$
et un prolongement analytique $(V_k, F_k)$
de
$F^k$ tel que $F_k(V_k)\subset\R$. Par cons\'equent,
pour tout $k\geq 1$, 
$h^{\cIrc k}(\R)\subset BPA(F,0)$.
\labelle{fixattr}
\end{lemme}
{\bf D\'emonstration:}
$(\Rightarrow)$ \'etant trivial, on prouvera $(\Leftarrow)$.
Pour chaque $0<\nu\leq k$
et $x\in h^{\nu}(\R)$ il existe
un inverse locale $\phi_{\nu,x}$ de $h^{\nu}$
et un voisinage ${\cal U}_{\nu,x}$ de $x$ tels que
$\phi_{\nu,x}\left({\cal U}_{\nu,x}    \right)\subset
{\R} $. Posons $x_0:=p $,
$
x_{\nu+1}:=
h^{k-(\nu-1)}
\circ \phi_{k-\nu,x_{\nu}}(x_{\nu})$
et $F_{\nu}:=h^k-(\nu-1)\circ \phi_{k-\nu,x_{\nu}}$
 pour 
$0\leq\nu<k$.
Alors $F_{\nu}(x_{\nu})=x_{\nu+1}$,
$F_{\nu}({\cal U}_{\nu,x_{\nu}})\subset{\cal U}_{\nu-1,x_{\nu+1}} $
et $\CIRC_{\nu=0}^k F_{\nu}({\cal U}_{k,p})\subset\R$.
On conclut en posant,
selon la notation
de la d\'efinition \ref{basattr}, $V_{\nu}:={{\cal U}_{k-\nu,x_{\nu}}}$.
\QUAN
\vskip0,2truecm
Nous rappelons que un {\sf domaine de Fatou-Bieberbach} est un ouvert de  $\CI^N   $
biholomorphe \`a  $\CI^N   $, peut 
\^etre un sous-ensemble propre.
On va d\'efinir une g\'en\'eralisation
\`a  'plusieurs feuillets' de cette id\'ee:
\begin{defi}
Un {\sf domaine de Riemann-Fatou-Bieberbach}
au-dessus d'un ouvert   $\Omega\subset\CI^N$
est un domaine de Riemann   $R$ au-dessus
de $\Omega$, avec projection $\pi$ qui est \'egalement un domaine de Riemann au-dessus 
de $\CI^N$ avec projection $\Psi$, tel que
$\Psi(x)=\Psi(y)\Rightarrow \pi(x)=\pi(y)$.
\labelle{rifatbib}
\end{defi}

\begin{theoreme}
Soit $h   $  un endomorphisme de $\CI^N$
avec un point fixe r\'epulsif r\'egulier en $0   $,
$\R$ un voisinage de r\'egu;arit\'e de $0$,
$\Omega   $ le  bassin potentiel 
d'attraction de $0$ pour la dynamique de
$F:=[h\vert_{\R}]^{-1}$. Alors il existe
un domaine de Riemann-Fatou-Bieberbach $\M$ 
au-dessus de $\Omega$.
De plus,
 si $0$ est fortement r\'egulier, $\M$
est biholomorphe \`a $\CI^N$.
\labelle{rudplus}
\begin{remarque}
On verra plus tard (voir section {\bf \ref{esempi}}) que, g\'en\'eralement,
on a $\Omega\not=\CI^N$, comme dans le 
cas des domaines de Fatou-Bieberbach.
\labelle{fabi}
\end{remarque}
\end{theoreme}
{\bf D\'emonstration}:
rappelons que 
$F^n(\R)\subset\alpha^n
\R$ pour
$0<\alpha<1$ convenable.

Gr\^ace au th\'eor\`eme de Poincar\'e-Dulac, on prouve,
comme 
dans la
la d\'emonstration
du {\tt th\'eor\`eme} de l'ap\-pen\-di\-ce de \cite{rudin}, 
qu'il existe une une  application holomorphe
$\Ppsi_0:\R\rightarrow \CI^N$ 
satisfaisant
\begin{equation}
\cases
{
\Ppsi_0=0\cr
(\Ppsi_0)_*=\id\cr
G^{-n}\hbox{\tt o}\Ppsi_0=\Ppsi_0\hbox{\tt 
o}\,\h^{n}.
}
\labelle{ast8}
\end{equation}

Consid\'erons maintenant le prolongement
analytique maximal 
$\left(
\M,\pi,j
,\Psi    \right)   $ de $\Ppsi_0$.

\vskip0.2truecm

\begin{lemme}
Si $x_1,x_2\in\M$ et $\Psi(x_1)=\Psi(x_2)$,
on a 
$\pi(x_1)=\pi(x_2)$.
\labelle{presquin} 
\end{lemme}
{\bf D\'emonstration:}
comme $\pi(x_1),\pi(x_2)\in\Omega$, il
existe
des voisinages ${\cal U}_i$ de $\pi(x_i)$,
des points $\{x_{ik}\}_{k=0...N}$
 (avec $x_{i0}=\pi(x_i)$, $i=1,2   $)
et
des \'el\'ements d'applications holomorphes $(W_{ik},f_{ik})$, $i=1,2   $, 
chacun desquels reliable avec
$F$, 
tels que
$$
\cases{
x_{ik}
\in 
W_{ik}
\cr
f_{ik}(x_{ik})=x_{i,{k+1}}
\cr
\lim_{k\to\infty}
\hbox{\tt\large o}_{l=0}^{N}\,
f_{il} (W_i0)\subset\R
}
\quad 
i=1,2.
$$
En outre, on peut supposer, sans perte de g\'en\'eralit\'e, que
$j$ admet des prolongements analytiques
$j_{ik}$ sur tous les $W_{ik}$, 
de fa\c con telle que 
$\{j_{ik}(W_{ik})\}$ (i=1,2) soient deux cha\^\i nes d'ouverts en $\M$ connectant
respectivement $j(0)$ avec $x_1$ et 
$j(0)$ avec $x_2$.
On peut aussi supposer que, pour $k$ assez grand,
$W_{ik}\equiv\R$ et ${j_{ik}\equiv j}$, $i=1,2$.

Posons $F_{ik}=\hbox{\tt\large o}
_{l=0}^{k}\,f_{il}$, $i=1,2$; on a
$\Psi\circ j_{ik}\circ F_{ik}\circ\pi=
G^k\circ\Psi$, donc
\begin{equation}
\Psi j_{1k} F_{1k}
\pi(x_1)= G^{k}\Psi(x_1)
=G^{k}\Psi(x_2)
=\Psi\ j_{2k} F_{2k}\pi(x_2).
\labelle{quasifine}
\end{equation}

On peut supposer, sans perte de g\'en\'eralit\'e, que $j$ soit inversible dans
$\hbox{\tt S}$ 
et $\Psi$ injective dans $j(\hbox{\tt S})$.
On a alors que (\ref{quasifine}) entra\^\i ne, pour $k=N$
$
j_{1N}\circ F_{1N}\circ\pi(x_1)=
j_{2N}\circ F_{2N}\circ\pi(x_2)
$;
en appliquant $h^{\cIrc N}\circ\pi$ on obtient
$\pi(x_1)=\pi(x_2)$.
\QUAN{\scriptsize\ (lemme \ref{presquin})}

\vskip0,2truecm
\begin{lemme}
On a:  {\tt(i)} $\Omega=\pi(\M)$;
{\tt (ii)} $\Psi   $ est localement biholomorphe;
{\bf si le point fixe $0$ est fortement r\'egulier},
{\tt (iii)} $\pi$ est un rev\^etement topologique et
{\tt (iv)} $\psi$ est un rev\^etement topologique. 

\labelle{dedans}
\end{lemme}
{\bf D\'emonstration.}
{\tt(i)}: gr\^ace au lemme \ref{fixattr}, $\Omega=\bigcup_{k=1}^{\infty}
\left[h^{\cIrc k}(\hbox{\R})    \right]$.
Comme $\R\subset h(\R)$
on a $h(\Omega)=\Omega$.
Prouvons d'abord que $\Omega\subset\pi(\M)$.
Si $p\in\Omega$, il existe, gr\^ace au lemme \ref{basattr},
$n\in\ENNE$ et un prolongement analytique $(V,F_n)$ de $(\R,F^n)$ \`a un voisinage
 $V$ de $p$, tel que $F_n(V)\subset\R$.
 
Cela entra\^\i ne que
la suite d'applications holomorphes
\begin{equation}
\Big\{
(G^{-k}\circ T\circ F^{\cIrc k-n})\hbox{\tt o}
F_n
\Big\}_{k\in\ENNE}
\labelle{compact}
\end{equation}
converge
uniform\'ement
sur 
les compacts de $V$ vers
une application holomorphe
$\Ppsi_p$, qui est visiblement un prolongement analytique de
$\Psi_0  $  car  
$F_n$ est un prolongement analytique de $F^n$. Ainsi $p\in\pi(\M)$.

Soit $p\in\pi(\M)$: il existe un prolongement analytique 
$({\cal U}_p, \Psi_p)$ de $\Ppsi_0$ \`a un voisinage ${\cal U}_p$ de $p$.

On peut supposer $\Ppsi_0$
inversible sur un ouvert ${\cal U}\subset\R$:
en soit $({\cal V},\Psi_0^{-1})$
l'invers.
On peut toujours supposer 
${\cal V}\subset\R$.
Comme 
\begin{equation}
\Psi_0^{-1}\circ G^n\circ\Psi_0=F^n
\label{local}
\end{equation}
sur $\R$
pour
tout $n$ et, par le lemme \ref{erreuno},
$\lim_{k\to\infty}G^{\cIrc k}=0$
uniform\'ement sur les compacts 
de $\CI^N$, 
on voit que, pour
chaque compact ${\cal K}\subset {\cal U}_p$ 
et $n$ assez grand, $G^n \Psi_p({\cal K})\subset\R$.
On peut  donc 
prolonger le membre gauche de (\ref{local}) sur ${\cal U}_p$, en 
gagnant l'\'el\'ement $\Psi_0^{-1}\circ G^n\circ\Psi_p$. Par cons\'equent, 
$F^n$ aussi peut \^etre prolong\'e sur ${\cal U}_p$ \`a un \'el\'ement 
$({\cal U}_p,F_n)$
et $F_n(p)\in\R$.
Ainsi $p\in h^{\cIrc n}(\R)$ et,
gr\^ace au lemme \ref{fixattr}, $p\in\Omega$.

{\tt (ii)}: 
notons que 
la d\'efinition de $\Ppsi_p$ par
(\ref{compact}) entra\^\i ne 
que cette application est une limite de biholomorphismes locaux, donc soit 
$\Ppsi_p$ est d\'eg\'en\'er\'ee au voisinage de $p$, soit elle y est biholomorphe.
Le premier cas ne peut pas se pre\'esenter, car sinon, par prolongement analytique,
m\^eme $\Ppsi_0$ serait  d\'eg\'en\'er\'ee, 
ce qui
contredit
(\ref{ast8}).

{\tt (iii)}: prouvons que, dans le cas de forte r\'egularit\'e, $\pi$ jouit de la propri\'et\'e du
rel\`evement des courbes.
Soit $\gamma:I\rightarrow\pi(\M)$ un chemin et $x\in\pi^{-1}(\gamma(0))$.
Par construction de $\M$ il existe un chemin $\beta:I\rightarrow\pi(\M)$ tel que
$\beta(0)=0$, $\beta(1)=\gamma(0)$ admettant un rel\`evement $\tilde\beta:I\rightarrow\M$
tel que $\tilde\beta(0)=j(0)$ et $\tilde\beta(1)=x$.
Soit $\Gamma:=\beta*\gamma$: comme $\pi(\M)=\Omega$ par {\tt (i)}, $\Omega=\bigcup_{k=1}^{\infty}
\left[h^{\cIrc k}(\hbox{\R})    \right]$ et
$h^k(\R)\subset h^{k+1}(\R)$, on a $\Gamma(I)\subset h^N(\R)$
pour $N$ assez grand.

Or $h^N$ est un rev\^etement, donc $(\R,F^N)$ admet un prolongement analytique
le long de $\Gamma$ jusqu'\`a un \'el\'ement $(V,F_n)$ dans un voisinage de $\gamma(1)=\Gamma(1)$.

On peut d\'efinir un prolongement analytique de $\Ppsi_0$ le long de $\Gamma$ \`a l'aide de
(\ref{compact}).

Cela entra\^\i ne que $\Gamma$ admet un rel\`evement
$\widetilde\Gamma:I\rightarrow(\M)$ tel que $\widetilde\Gamma(0)=j(0)$ et 
$\widetilde\Gamma(1/2)=x$.
Posons $\tilde\gamma(t):=\widetilde\Gamma((t+1)/2)$: on voit que $\tilde\gamma(0)=x$ et
$\pi\tilde\gamma=\gamma$, donc $\gamma$
admet un rel\`evement respectivement \`a $\pi$ commen\c cant \`a $x$, c'est-\`a-dire $\pi$
est un rev\^etement topologique.

{\tt (iv)}: au voisinage de $0$ on a 
\begin{equation}
\Ppsi_0^{-1}=\lim_{k\to\infty} h^k\circ T^{-1}\circ G^k;
\labelle{duesestelles}
\end{equation}
cette d\'efinition-l\`a peut \^etre prolong\'ee \`a une application holomorphe $\Theta$
sur $\CI^N$, car $\lim_{k\to\infty}G^k=0$ uniformement sur les compacts de $\CI^N$.

Or, au voisinage de chaque point $p\in\CI^N$, la suite (\ref{duesestelles})
est une suite de biholomorphismes locaux, car les $\{h^k\}$ le sont 
sur $\R$, donc soit $\Theta$ est d\'eg\'en\'er\'ee au voisinage de $p$ soit $p$ n'est pas un point critique
pour $\Theta$.
Le premier cas ne peut pas se pre\'esenter, car sinon $\Theta$ serait d\'eg\'en\'er\'ee sur $\CI^N$, ce qui contredit $\Theta_*(0)=\Ppsi_*(0)^{-1}=\id$.
Gr\^ace au lemme \ref{inverse} et \`a {\tt (ii)}, $\Theta(\CI^N)=\Omega$.

Montrons que $\psi$ jouit de la propri\'et\'e du
rel\`evement des courbes.
Soit $\gamma:I\rightarrow\CI^N$ un chemin et $y\in\psi^{-1}(\gamma(0))$: gr\^ace au lemme 
\ref{presquin}, $y\in\pi^{-1}(\Theta(\gamma(0)))$.
Puisque $\pi$ est, par {\tt (iii)}, un rev\^etement, il existe un rel\`evement
$\tilde\gamma:I\rightarrow\M$
de $\Theta\gamma$ commen\c cant \`a $y$; comme on a aussi $\psi\tilde\gamma=\gamma$, on voit que 
$\tilde\gamma$ est un rel\`evement de $\gamma$ par respectivement  \`a $\Psi$, commen\c cant \`a $y$, c'est \`a
dire $\Psi$ est un rev\^etement topologique.

\QUAN{\scriptsize\ (lemme \ref{dedans})}

\vskip0.2truecm

{\bf Fin de la d\'emonstration du th\'eor\`eme \ref{rudplus}:}
montrons que $\Psi   $
est surjective:
on peut recouvrir $\M$ par un ensemble
d\'enombrable
 d'ouverts
$\{{\cal V}_l\}$ tels que $\pi\vert_{{\cal V}_l}$
est inversible; posons 
$
{\cal U}_l
:=
\pi({{\cal V}_l})
$.
Gr\^ace au lemme \ref{dedans} {\tt (i)} et au fait que $h(\Omega)=\Omega$,
pour tout $n$ et tout $l$ on aura aussi
$h^n({\cal U}_l)\subset\bigcup_{\lambda\in
L(l,n)}{\cal U}_l$, pour un certain ensemble d'indices $L(l,n)$.
Par construction $\Psi(M)=
\bigcup_{l=1}^{\infty}\Psi_l
({\cal U}_l)$, o\`u chacun des 
$({\cal U}_l,\Psi_l)$ est un prolongement analytique de
$\Ppsi_0$, et donc
$$
G^{-n}\Psi(\M)=
\bigcup_{l=1}^{\infty}
G^{-n}
\Psi_l
({\cal U}_l)
.
$$

Par prolongement analytique de (\ref{ast8}), on a, 
pour chaque $l,n$,
$$
G^{-n}
\Psi_l
({\cal U}_l)
=
\bigcup_{\lambda\in L(l,n)}
\Psi_{\lambda}
h({\cal U}_l)
\subset
\Psi(M);
$$
en consid\'erant la r\'eunion sur $l$, 
on obtient
$G^{-n}\Psi(M)\subset\Psi(M)$; mais alors, gr\^ace au lemme \ref{erretre},
$$
\CI^N=\bigcup_{n=1}^{\infty}
G^{-n}\Psi(M)\subset\Psi(M),
$$
donc 
$\Psi(M)=\CI^N$.

Finalement, {\bf dans le cas fortement r\'egulier}, gr\^ace au lemme \ref{dedans} {\tt (iv)},
$\Psi$ est un rev\^etement topologique, donc, gr\^ace \`a la connexit\'e simple de $\CI^N$, 
il s'agit d'une application biholomorphe.

\QUAN{\scriptsize \ (th\'eor\`eme \ref{rudplus})}
\section{Exemples \labelle{esempi}}
\subsection{Le cas r\'egulier}

Voici une classe 
d'exemples de 
domaine de Riemann-Fatou-Bieberbach:
on construira d'abord 
un automorphisme $h$ de $\CI^N$, 
admettant 
un point fixe r\'epulsif r\'egulier en $0   $,
et
dont 
l'image soit un sous-ensemble ouvert
propre $\Omega$ de $\CI^N$.
Pour ce faire,
il suffit de prendre en consid\'eration
une application biholomorphe
$\chi:\CI^N\rightarrow \CI^N$, tangente a l'identit\'e
en $0$, 
telle que $\chi(\CI^N)\,(\not=\CI^N)$
soit 
le bassin d'attraction
du point fixe attractif $0$ pour 
un automorphisme
de $\CI^N$.
En posant
$h=2\chi$ et $S:=\chi(\CI^N)$, nous
obtenons une application biholomorphe
$h:\CI^N\rightarrow 2S$, dont
nous appellerons $f$ l'inverse.
Comme $f$ a un point fixe attractif en $0$,
il existe un bassin d'attraction
$\Omega$ de ce point.
On a forc\'ement $\Omega\subset
2S$.

Soit maintenant $N=2$ et ${\cal E}:
\CI^2\rightarrow\CI^2$
d\'efinie en posant 
${\cal E}(u,v)=
(\exp(u)-1,v)$ et
$\Ef:={\cal E}\circ h.$
On a
$\Ef(0)=0$ et $\Ef_*(0)=h_*(0)=2\id$.
On voit ais\'ement que $0$ est r\'epulsif r\'egulier
pour $\Ef$.
Si $\psi$ est une inverse locale de $\Ef$ en $0$,
on a $\vert\vert \psi_*(0)\vert\vert=1/2$,
donc on peut it\'erer 
$\psi$ et il existe un voisinage ouvert 
${\cal U}$ de
$0$ tel que 
$
K\subset\subset {\cal U}\Rightarrow 
\psi^{\cIrc n}(K)\rightarrow 0
$,
c'est-\`a-dire, ${\cal U}$ est dans le bassin
potentiel d'attraction
$\Omega^{\prime}$
de $0$ pour la dynamique de $\psi$.

Gr\^ace au lemme \ref{fixattr},
$\Omega^{\prime}=\bigcup_{n=1}^{\infty}
\Ef^{k}({\cal U})
$; comme 
$
{\cal E}(2S)=\phi(\CI^2)
\supset \phi^{\cIrc 2}(\CI^2)\supset
\cdots\supset \phi^{\cIrc n}(\CI^2)
\supset\cdots$,
on a 
$
\Omega^{\prime}\subset
\subset{\cal E}(2S):
$
c'est un sous-ensemble propre de $\CI^2$.
\subsection{Le cas 
fortement 
r\'egulier}
Soit ${\cal E}:=\CI^2_{(z,w)}
\setminus 
\{0\}\times\CI$: 
construissons
une application
holomorphe propre 
$h:\CI^2\rightarrow{\cal E}
$ 
avec un point fixe r\'epulsif
et
$Jh\not=0$ sur 
${\cal E}$.
Pour ce faire
on pourra
partir de l'application biholomorphe
sur $\CI^2$, prenant valeurs
en
${\cal E}$, de \cite{rudin}, p.76
o\`u 77, que nous 
allons appeler
$G$.
Soit $G(1,1)=(a,b)
(\in{\cal E})$;
or, pour $\alpha$, $\beta$
convenables l'application
$H$ d\'efinie en posant
$$
H(z,w):= G(\alpha z-
\alpha a+1, \beta y -
\beta b+1)
$$
a
un point fixe
r\'epulsif en $p:=(a,b)$.
Enfin, on pourra  
considerer 
$h:=H(a^{1-n} z^n,w)$, qui
jouit 
\'evidemment
des propri\'et\'es
\'enonc\'ees au debut du 
paragraphe.

Or, pour chaque voisinage $V$
de $p$ en ${\cal E}$
et pour chaque entier
positif $k$, on a 
$h^k(V)\subset{\cal E} $ 
par construction.
Par ailleurs, $h^k\vert_{V}$
est un biholomorphisme
local propre, et donc un 
rev\^etement topologique.
Ainsi le point fixe $p$
est fortement r\'egulier.
Donc la construction
du th\'eor\`eme \ref{rudplus}
nous donne un domaine de Riemann
 $\M$ biholomorphe \`a $\CI^N$
qui recouvre le bassin potentiel
d'attraction de $p$ par respect
\`a
un inverse local de $h$.
Par construction, ce 
bassin est contenu
en ${\cal E}$, et donc il est
un sousensemble propre de
$\CI^2$.

\end{document}